\documentclass{amsart}

\usepackage{amsmath,amssymb,amsthm}

\usepackage{graphicx,tikz}

\newtheorem*{thm}{Theorem}

\theoremstyle{definition}

\theoremstyle{remark}

\begin{document}

\title[]{An Agmon estimate for \\Schr\"odinger operators on Graphs}
\thanks{S.S. is supported by the NSF (DMS-2123224) and the Alfred P. Sloan Foundation.}

\author[]{Stefan Steinerberger}
\address{Department of Mathematics, University of Washington, Seattle, WA 98195, USA}
\email{steinerb@uw.edu}

\begin{abstract} The Agmon estimate shows that eigenfunctions of Schr\"odinger operators,
$ -\Delta \phi + V \phi = E \phi$, decay exponentially in the `classically forbidden' region where the potential exceeds
the energy level $\left\{x: V(x) > E \right\}$. Moreover, the size of $|\phi(x)|$ is bounded in terms
of a weighted (Agmon) distance between $x$ and the allowed region. We derive such a statement on graphs when $-\Delta$ is replaced by the Graph Laplacian $L = D-A$: we identify an explicit Agmon metric and prove a pointwise decay estimate in terms of the Agmon distance. 
 \end{abstract}

\maketitle

\section{Introduction: Agmon estimates}
Let $V: \mathbb{R}^n \rightarrow \mathbb{R}$ be a non-negative potential growing at infinity, i.e. $V(x) \rightarrow \infty$ as $\|x\| \rightarrow \infty$. 
Agmon estimates are concerned with eigenfunctions of the Schr\"odinger operator $ -\Delta + V$: we study functions $\phi: \mathbb{R}^n \rightarrow \mathbb{R}$ satisfying
$$ - \Delta \phi + V \phi = E \phi,$$
where $E \in \mathbb{R}$ is the eigenvalue. Multiplying with $\phi$ and integrating by parts,
$$ \int_{\mathbb{R}^n} |\nabla \phi|^2 dx + \int_{\mathbb{R}^n} V(x) \cdot \phi(x)^2 dx =  \int_{\mathbb{R}^n} E \cdot \phi(x)^2 dx.$$
This identity implies that most of the $L^2-$mass should be contained in the `allowed' region $\left\{x \in \mathbb{R}^n:  V(x) \leq E \right\}$ and
only very little mass can be in the `forbidden' region $\left\{x \in \mathbb{R}^n: V(x) > E \right\}$. 
The celebrated Agmon estimate \cite{agmon1} shows that this is indeed the case and, moreover, that $\phi$ decays exponentially in terms of the distance from the allowed region for a suitable notion of distance. 
Agmon's estimate can be derived from an explicit integral identity. One way of motivating the estimate (taken from a summary of Deift \cite{deift}) is as follows: if
$$ - \Delta \phi + V\phi = E \phi,$$
then for any (sufficiently regular) $g: \mathbb{R}^n \rightarrow \mathbb{R}$ 
$$ \int_{\mathbb{R}^n} \left|  \nabla  (e^{ g} \phi) \right|^2 dx + \int_{\mathbb{R}^n} \left(V - E - | \nabla g|^2 \right) e^{2g} \phi^2 dx = 0.$$
Ignoring the first (positive) term, this implies
$$  \int_{\mathbb{R}^n} \left(V - E - | \nabla g|^2 \right) e^{2g} \phi^2 dx \leq 0.$$
We note that $e^{2g}$ and $\phi^2$ are positive, $V - E - | \nabla g|^2$ is negative in the allowed region and positive in the forbidden region provided $|\nabla g|$ is sufficiently small. The inequality then naturally implies that there cannot be too much $L^2-$mass of $\phi$ in the forbidden region except this is now coupled with an additional exponentially growing term $e^{2g}$. The statement becomes stronger the larger we make $g$, however, we want to maintain the nonnegativity of $V-E- |\nabla g|^2$ in the forbidden region. This then suggests a way of defining $g$: the \textit{Agmon metric} associated to the energy level $E$ between two points $x,y \in \mathbb{R}^n$ is given as the minimum
energy taken
$$ \rho_{E}(x,y) = \inf_{\gamma} \int_0^1 \max\left( \sqrt{V(\gamma(t))- E}, 0\right) |\dot \gamma(t)| dt,$$
where $\gamma:[0,1] \rightarrow \mathbb{R}^n$ ranges over all paths from $\gamma(0) = x$ to $\gamma(1) = y$. The integral identity can then be used (see for example Carmona \& Simon \cite{carmona}) to derive pointwise statements in the forbidden region along the lines of that for all $\varepsilon >0$
$$ |\phi(x)| \leq c_{\varepsilon} \sup_{y \in \mathbb{R}^n \atop V(y) \leq E} e^{-(1-\varepsilon) \rho_{E}(x,y)}.$$
We refer to Aizenman \& Simon \cite{aiz}, Carmona \cite{carmona0}, Dimassi \& Sj\"ostrand \cite{dim}, Helffer \cite{helffer}, Helffer \& Sj\"ostrand \cite{helff2, helff3}, Hislop \cite{hislop}, Simon \cite{simon, simon3} and references therein for a more complete picture regarding Agmon's estimate in
the continuous setting. Our paper is partially inspired by a recent probabilistic approach to obtain sharp pointwise Agmon estimates in the continuous setting \cite{steini}.

\section{An Agmon estimate on graphs}
\subsection{Setup.}
Let $G = (V,E)$ be a finite, connected graph with $V = \left\{v_1, \dots, v_n\right\}$. We introduce the diagonal matrix $D \in \mathbb{R}^{n \times n}$ satisfying $d_{ii} = \deg(v_i)$ and the adjacency matrix $A \in \mathbb{R}^{n \times n}$ given by
$$ A_{ij} = \begin{cases} 1 \qquad &\mbox{if}~(v_i, v_j) \in E \\ 0 \qquad &\mbox{otherwise.} \end{cases}$$
There is a natural notion of a discrete Laplacian acting on functions $f:V \rightarrow V$  given by the linear operator $L= D-A \in \mathbb{R}^{n \times n}$. We observe that $L$ can be interpreted as the discrete analogue of $-\Delta$ as both are positive semi-definite and allow for integration by parts: for $f:V \rightarrow \mathbb{R}$
$$  \left\langle f, L f \right\rangle = \sum_{(v_i, v_j) \in E} (f(v_i) - f(v_j))^2.$$
Given an arbitrary potential $W: V \rightarrow \mathbb{R}_{}$,
our goal is to understand the behavior of eigenfunctions $\phi:V \rightarrow \mathbb{R}$ satisfying
$$ L \phi + W \phi = E \phi$$
for some eigenvalue $E \in \mathbb{R}$.
Multiplying with $\phi$ and integrating by parts
$$ \sum_{(v_i, v_j) \in E} (\phi(v_i) - \phi(v_j))^2 + \sum_{v \in V} W(v) \cdot \phi(v)^2 =  \sum_{v \in V} E \cdot \phi(v)^2.$$
This suggests, just as in the continuous case above, that there should be relatively little $\ell^2-$mass in the `classically forbidden' region
$\left\{v \in V: W(v) > E \right\}$. The question is now whether, just as in the continuous case, one can expect exponential decay in the forbidden region and how this can be quantified.

\subsection{Main Result.} We define a notion of Agmon distance $\rho_{E}:V \rightarrow \mathbb{R}$ as the cost of the cheapest path starting in $v \in V$ and ending in any arbitrary vertex in the allowed region where `cheap' refers to an explicit cost function on $V$ depending on the potential $W$, the energy $E$ and the degree of the vertex. Formally, 
$$ \rho_E(v) = \inf \left\{\sum_{i=1}^{\ell}  \log\left(1 +  \frac{(W(v_i)-E)_{+}}{\deg(v_i)} \right):  v = v_1 \rightarrow \dots \rightarrow v_{\ell} ~ \mbox{and}~~ W(v_{\ell}) \leq E \right\},$$
where the infimum is taken over all paths that start in $v$ and end in a vertex $v_{\ell}$ in the allowed region. As usual, $(W(v_i) - E)_{+} = \max\left\{ W(v_i) - E, 0 \right\}$. 
Note that $\rho_E \equiv 0$ in the allowed region and $\rho_E > 0$ in the forbidden region. 
\begin{thm}
We have, for all $v\in V$,
$$ |\phi(v)| \leq e^{ - \rho_E(v) } \cdot \| \phi\|_{\ell^{\infty}}.$$
\end{thm}


$|\phi|$ assumes its maximum in the allowed region and thus
$ \| \phi\|_{\ell^{\infty}} = \| \phi\|_{\ell^{\infty}(W(v) \leq E)}.$
Since $\rho_E \equiv 0$ in the allowed region, the inequality is sharp in the maximum and the implicit constant 1 in front cannot be improved any further. In terms of the exponential decay,
there are graphs where the inequality is asymptotically optimal: such examples are constructed in \S 3.2.

\subsection{Related results.}
There is relatively little work regarding Agmon estimates on graphs. However, we emphasize one recent result which is close in spirit to our result.  Filoche, Mayboroda \& Tao \cite{filoche} study eigenvector localization for a fairly general class of matrices $A \in \mathbb{R}^{n \times n}$. They obtain an integrated exponential estimate in terms of an explicit Agmon-type distance. Considering $A = D-A +W$ and $u = (1,1,\dots,1)$ in their approach, one arrives at a notion of distance
 $$ \rho(v, w) = \inf \left\{\sum_{i=1}^{\ell} \log\left(1 + \sqrt[4]{(W(v_i) -E)_{+}(W(v_{i+1})-E)_{+}}\right) \right\},$$
 where the infimum ranges over all paths that start in $v=v_1$ and end in $w=v_{\ell+1}$ (and, as in our approach, traveling through the allowed region is free which we suppress in the equation above for simplicity of exposition). This is very similar in flavor to our distance above: using this, they then
 obtain an integrated estimate also involving a landscape-type potential $A^{-1} \mathbf{1}$ \cite[Theorem 2.5]{filoche} as well as more general integrated estimates \cite[Theorem 2.7]{filoche}. A main
 difference is the dependency on the degree of a vertex which is locally built into our distance while arising in the integrated estimates of \cite{filoche} more globally (somewhat unsurprisingly: integrated estimates themselves are global). Both our estimate and the estimates in \cite{filoche} are complementary: which one ends up being better will depend (among other things) on whether there is a lot of variation in the degrees of the vertices.\\
 
There is also a recent work of Keller \& Pogorzelski \cite{keller} who study Agmon estimates in the more general setting of weighted, infinite graphs where the Agmon distance is given in terms of Hardy weights.
We also note the work of Akduman \& Pankov \cite{quantum, quantum2} on metric graphs, the work of Harrell \& Maltsev \cite{evans} on quantum graphs, 
Damanik, Fillman \& Sukhtaiev \cite{damanik} on tree graphs
and Klein \& Rosenberger \cite{klein, klein2}, Mandich \cite{man} and  Wang \& Zhang \cite{wang}  on $\mathbb{Z}^d$.

\section{Proof}
\subsection{Proof of the Theorem}  
%
%

\begin{proof}
Note first that
the eigenfunction satisfies
$$ (D-A)\phi = (E-W) \phi.$$
Considering this linear system of equations in a fixed vertex $u \in V$ one obtains
$$ \deg(u) \phi(u) - \sum_{(u,w) \in E} \phi(w) = (E-W(u)) \phi(u).$$
This equation can be rewritten as
$$ \left[ 1 + \frac{W(u)-E}{\deg(u)} \right] \phi(u) = \frac{1}{\deg(u)}\sum_{(u,w) \in E} \phi(w).$$
We observe that if $\phi(u) = 0$, then the Theorem is trivially true in $u$. It thus suffices to prove it for
vertices $u \in V$ where $\phi(u) \neq 0$. Note, moreover, that in the forbidden region $\left\{u \in V: W(u) > E\right\}$, one trivially has 
$$ 1 + \frac{W(u)-E}{\deg(u)} \geq 1$$
and thus, for $u \in V$ in the forbidden region, it is possible to divide and
$$ \phi(u) = \left[ 1 + \frac{W(u)-E}{\deg(u)} \right]^{-1}  \frac{1}{\deg(u)}\sum_{(u,w) \in E} \phi(w).$$
Taking absolute values on both sides, we have
\begin{align*}
 |\phi(u)| &\leq  \left[ 1 + \frac{W(u)-E}{\deg(u)} \right]^{-1}  \frac{1}{\deg(u)}\sum_{(u,w) \in E} |\phi(w)| \\
 &\leq \left[ 1 + \frac{W(u)-E}{\deg(u)} \right]^{-1} \max_{(u,w) \in E} |\phi(w)|.
 \end{align*}
Since $\phi(u) \neq 0$, we deduce 
$$  \max_{(u,w) \in E} |\phi(w)| > |\phi(u)|.$$
We can now move from $u$ to its neighbor $w$ maximizing $|\phi(w)|$ and then apply the very same argument again in $w$. The argument can be applied iteratively 
 as long as the new vertex is still in the forbidden region. Note that $|\phi|$ is increasing along the way which implies that the arising path can never cross itself and must eventually
end up in the allowed region $\left\{v \in V: W(v) \leq E\right\}$. Altogether, this results in a path
$$ u = v_1 \rightarrow v_2 \rightarrow \dots \rightarrow v_m \rightarrow v_{m+1} \qquad \mbox{where} \quad W(v_{m+1}) \leq E$$
since otherwise the path could be further extended.
Collecting all the factors 
\begin{align*}
 |\phi(u)| &\leq \left( \prod_{i=1}^{m}   \left[ 1 + \frac{W(v_i)-E}{\deg(v_i)} \right]^{-1} \right) \cdot |\phi(v_{m+1})| \\
  &\leq \left( \prod_{i=1}^{m}   \left[ 1 + \frac{W(v_i)-E}{\deg(v_i)} \right]^{-1} \right) \cdot \| \phi\|_{\ell^{\infty}}
 \end{align*}
Note that
$$   \prod_{i=1}^{m}   \left[ 1 + \frac{W(v_i)-E}{\deg(v_i)} \right]^{-1} = \exp\left( - \sum_{i=1}^{m} \log\left(1 +  \frac{W(v_i)-E}{\deg(v_i)} \right) \right).$$
By definition of $\rho_E$, we have
$$  \exp\left( - \sum_{i=1}^{m} \log\left(1 +  \frac{W(v_i)-E}{\deg(v_i)} \right) \right) \leq \exp\left( - \rho_E(v_1)\right)$$
and this concludes the proof.
\end{proof}

\textbf{Remark.} We note that the final estimate in the argument implies 
$$  |\phi(u)| \leq  \exp\left( - \sum_{i=1}^{m} \log\left(1 +  \frac{W(v_i)-E}{\deg(v_i)} \right) \right) \cdot |\phi(w_{m+1})|$$
for any path starting in $u = v_1$ and ending in the vertex $v_{m+1}$ in the allowed region. This would imply a slightly refined
estimate where one is not only interested in minimizing the Agmon metric but also wants to end up in a vertex in the allowed
region such that $|\phi(w_{m+1})|$ is not too small.\\

\textbf{Remark.} We quickly note a part in the derivation where the argument can be lossy: the main inequality is
\begin{align*}
 |\phi(u)| &\leq  \left[ 1 + \frac{W(u)-E}{\deg(u)} \right]^{-1}  \frac{1}{\deg(u)}\sum_{w \sim u} |\phi(w)| \\
 &\leq \left[ 1 + \frac{W(u)-E}{\deg(u)} \right]^{-1} \max_{(u,w) \in E} |\phi(w)|.
 \end{align*}
This inequality could also be interpreted as a type of martingale inequality
and can be exploited in this sense.
Let us define a sequence of random vertices given by $X_0 = u$ and such that $X_{k+1}$ is a randomly chosen neighbor
of $X_k$ and that this random walk is continued until $W(X_k) \leq E$. We will denote the smallest such $k$ by the stopping time $\tau$. Assume furthermore that
$$  \frac{W(u)-E}{\deg(u)} \geq \delta \qquad \mbox{for all}~u \in V~\mbox{in the forbidden region.}$$
An iterative application of the inequality then implies
$$ |\phi(u)| \leq \left(\sum_{\ell =0}^{\infty} \frac{\mathbb{P}\left(\tau = \ell\right)}{(1+\delta)^{\ell}}\right) \cdot  \|\phi\|_{\ell^{\infty}} .$$
We note that this inequality can lead to improved results in settings where a random walk needs a very long time before arriving in the allowed region.
Observe that the sum can be interpreted as an exponential moment $\mathbb{E} \exp\left(\tau/(\delta+1)\right)$ of the stopping time $\tau$ which is a well studied object. We refer to \cite{steini}
for the derivation of Agmon estimates via this more stochastic perspective in the continuous setting.

\subsection{An Example.}  The purpose of this section is to construct a graph where the inequality is nearly sharp. One example of such graphs is given by
$q-$regular trees of a certain depth where the final layer of vertices is then additionally connected
to another vertex $v_{*}$ (see Fig. 2 for an example). We consider the potential given by $W(v_*) = 0$ and, for all other vertices $v \neq v_*$, we choose the potential to be constant and
$W(v) = W \gg q^k \gg 1$ for some very large constant $W \in \mathbb{R}$. The function we will consider is the first eigenfunction of $L +W$.

\begin{center}
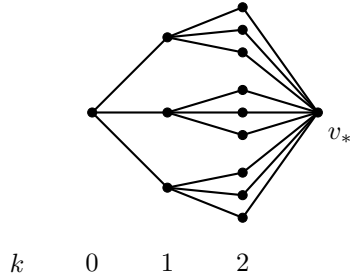
\begin{figure}[h!]
\begin{tikzpicture}
\filldraw (0,0) circle (0.06cm);
\filldraw (1,-1) circle (0.06cm);
\filldraw (1,0) circle (0.06cm);
\filldraw (1,1) circle (0.06cm);
\draw [thick] (0,0) -- (1,-1);
\draw [thick] (0,0) -- (1,0);
\draw [thick] (0,0) -- (1,1);
\filldraw (2,-1.4) circle (0.06cm);
\filldraw (2,-1.1) circle (0.06cm);
\filldraw (2,-0.8) circle (0.06cm);
\draw [thick] (1, -1) --  (2,-1.4);
\draw [thick] (1, -1) --  (2,-1.1);
\draw [thick] (1, -1) --  (2,-0.8);
\filldraw (2,-0.3) circle (0.06cm);
\filldraw (2,0) circle (0.06cm);
\filldraw (2,0.3) circle (0.06cm);
\draw [thick] (1, 0) --  (2,-0.3);
\draw [thick] (1, 0) --  (2,-0);
\draw [thick] (1, 0) --  (2,0.3);
\filldraw (2,1.4) circle (0.06cm);
\filldraw (2,1.1) circle (0.06cm);
\filldraw (2,0.8) circle (0.06cm);
\draw [thick] (1, 1) --  (2,1.4);
\draw [thick] (1, 1) --  (2,1.1);
\draw [thick] (1, 1) --  (2,0.8);
\filldraw (3,0) circle (0.06cm);
\draw [thick] (3, 0) --  (2,1.4);
\draw [thick] (3, 0) --  (2,1.1);
\draw [thick] (3, 0) --  (2,0.8);
\draw [thick] (3, 0) --  (2,-0.3);
\draw [thick] (3, 0) --  (2,-0);
\draw [thick] (3, 0) --  (2,0.3);
\draw [thick] (3, 0) --  (2,-1.4);
\draw [thick] (3, 0) --  (2,-1.1);
\draw [thick] (3, 0) --  (2,-0.8);
\node at (3.3, -0.3) {$v_*$};
\node at (-1, -2) {$k$};
\node at (0, -2) {$0$};
\node at (1, -2) {$1$};
\node at (2, -2) {$2$};
\node at (5,0) {};
\end{tikzpicture}
\caption{A $q-$regular tree (here, $q=3$) of depth $k$ (here, $k=2$) with final layer being connected to a single additional vertex $v_*$.}
\end{figure}
\end{center}

By Rayleigh-Ritz, the smallest eigenvalue of $-\Delta + W$
satisfies
$$ \lambda_1 = \inf_{f:V \rightarrow \mathbb{R}} \frac{ \sum_{(v_i, v_j) \in E} (f(v_i) - f(v_j))^2 + \sum_{v \in V} W(v) f(v)^2}{\sum_{v\in V}f(v)^2}.$$
Taking $f:V \rightarrow \mathbb{R}$ given by $f(v_*) = 1$ and $f(v) = 0$ for all $v \neq v^*$, we deduce that
$$ E = \lambda_1 \leq  q^k \qquad \mbox{independently of}~W.$$
We can now use the equation
$$ \phi(u) = \frac{E - W(u)}{\deg(u)} \phi(u) +  \frac{1}{\deg(u)}\sum_{w \sim u} \phi(w) $$
and we shall restrict its use to vertices in the $q-$regular tree. The value of $\phi$ then only depends on the level. We shall therefore write
$\phi(v) = \phi_i$ whenever the vertex $v$ is in the $i-$th level where $1 \leq i \leq k-1$ (the case $i=0$ and $i=k$ will be ignored since the algebra
is slightly different). The equation then simplifies to
$$ \phi_i = \frac{\lambda_1 - W}{q+1} \phi_i + \frac{1}{q+1}(q \cdot \phi_{i+1} + \phi_{i-1}).$$
This can be rewritten as
$$ (W - \lambda_1 + q + 1) \cdot \phi_i = q \cdot \phi_{i+1} + \phi_{i-1}.$$
For fixed $q$ and $W \gg q^k \geq E$, this implies that approximately $\phi_{i} \sim (q/W)\cdot \phi_{i+1}$ as $W \rightarrow \infty$
which implies exponential decay. Conversely, we have 
$$  \log\left(1 + \frac{W - E}{\deg(v)}\right) \sim \log\left(\frac{W}{q}\right)$$
which implies, to leading order, the same kind of decay.


\begin{thebibliography}{10}


\bibitem{agmon1} S. Agmon, Lectures on exponential decay of solutions of second order elliptic equations. Bounds
on eigenfunctions of N-body Schr\"odinger operators. Mathematical Notes, Princeton Univ. Press,
Princeton, N.J., 1982.

\bibitem{aiz} M. Aizenman and B. Simon, Brownian motion and Harnack's inequality for Schr\"odinger
operators, Comm. Pure Appl. Math. 35 (1982), 209--271.


\bibitem{quantum} S. Akduman and A. Pankov, 
Exponential estimates for quantum graphs. 
Electron. J. Differential Equations 2018, Paper No. 162, 12 pp.


\bibitem{quantum2} S. Akduman and A. Pankov, Nonlinear Schrödinger equation with growing potential on infinite metric graphs. Nonlinear analysis, 184 (2019),  p. 258--272.

\bibitem{carmona0} R. Carmona, 
Pointwise bounds for Schr\"odinger eigenstates.
Comm. Math. Phys. 62 (1978), no. 2, 97--106.

\bibitem{carmona} R. Carmona and B. Simon, Pointwise Bounds on Eigenfunctions
and Wave Packets in $N-$Body Quantum Systems, Commun. Math. Phys. 80, 59-98 (1981)

\bibitem{damanik} D. Damanik, J. Fillman \& S. Sukhtaiev, Localization for Anderson models on metric and discrete tree graphs. Mathematische Annalen, 376 (2020), p. 1337--1393.


\bibitem{deift} P. Deift, Review of \cite{agmon1}, Bull. Amer. Math. Soc. 12, p. 165--169 (1985) 


\bibitem{dim} M. Dimassi and J. Sj\"ostrand, Spectral asymptotics in the semi-classical limit, London Mathematical
Society Lecture Note Series, vol. 268, Cambridge University Press, 1999.

\bibitem{filoche} M. Filoche, S. Mayboroda and T. Tao,
The effective potential of an M-matrix.
J. Math. Phys. 62 (2021), no. 4, Paper No. 041902, 15 pp.

\bibitem{evans} E. Harrell and A.  Maltsev,
Localization and landscape functions on quantum graphs. 
Trans. Amer. Math. Soc. 373 (2020), no. 3, 1701--1729.

\bibitem{helffer} B. Helffer, Semi-classical analysis for the Schr¨odinger operator and applications.
Springer Lecture Notes in Math, 1988.

\bibitem{helff2} B. Helffer and J. Sj\"ostrand, Multiple wells in the semi-classical limit I, Comm. PDE 9, (1984), p. 337--408.


\bibitem{helff3} B. Helffer and J. Sj\"ostrand, Analyse semi-classique pour l’´equation de Harper (avec
application a l'equation de Schr\"odinger avec champ magn\'etique).
M\'emoire de la Soci\'et\'e Math\'ematique de France 34 (1988).

\bibitem{hislop} P. Hislop, Exponential decay of two-body eigenfunctions: A review, Proceedings of the Symposium on Mathematical Physics and Quantum Field Theory (Berkeley, CA, 1999). 

\bibitem{keller} M. Keller and F. Pogorzelski, Agmon estimates for Schr\"odinger operators on graphs, arXiv:2104.04737

\bibitem{klein} M. Klein and E. Rosenberger,
Agmon-type estimates for a class of difference operators. 
Ann. Henri Poincar\'e 9 (2008), no. 6, 1177--1215.

\bibitem{klein2} M. Klein and E. Rosenberger, The tunneling effect for a class of difference operators. 
Rev. Math. Phys. 30 (2018), no. 4, 1830002, 42 pp.


\bibitem{lithner} L. Lithner, 
A theorem of the Phragmen-Lindelof type for second-order elliptic operators.
Ark. Mat. 5 (1964), 281--285 (1964).

\bibitem{man} M-A. Mandich. Sub-exponential decay of eigenfunctions for some discrete Schr\"odinger operators. J. Spectr. Theory, 9(1):21--77, 2019.

\bibitem{simon} B. Simon. Semiclassical analysis of low lying eigenvalues, I. Nondegenerate minima: Asymptotic expansions,
Ann. Inst. H. Poincare 38 (1983), 295--307.

\bibitem{simon3} B. Simon, Semi-classical analysis of low lying eigenvalues II. Tunneling.
Ann. of Math. 120 (1984), p. 89--118.

\bibitem{steini} S. Steinerberger, Effective Bounds for the Decay of Schr\"odinger Eigenfunctions and Agmon bubbles, arXiv:2110.01163

\bibitem{wang} W. Wang and S. Zhang, 
The exponential decay of eigenfunctions for tight-binding Hamiltonians via landscape and dual landscape functions. 
Ann. Henri Poincar\'e 22 (2021), no. 5, 1429--1457.

\end{thebibliography}
\end{document}